\documentclass[]{amsart}%{article}
\usepackage{amssymb}
\newcommand{\Dl}{\Delta}
\newcommand{\Lm}{\Lambda}

\newcommand{\sg}{\sigma}

\newcommand{\lm}{\lambda}

\newcommand{\ra}{\rightarrow}

\newcommand{\sub}{\subset}
\newcommand{\nul}{\varnothing}

\newcommand{\Aut}{{\rm Aut}\,}

\newcommand{\N}{{\mathbb N}}%{{\bf F}}%
\newcommand{\R}{{\mathbb R}}%{{\bf R}}%
\newcommand{\T}{{\mathbb T}}%{{\bf T}}%
\newcommand{\Z}{{\mathbb Z}}

\newcommand{\tw}{\textwidth}

\newcommand{\bmp}{\begin{minipage}[t]{.9\tw}}
\newcommand{\emp}{\end{minipage}}

\theoremstyle{plain}%{definition}%{remark}%
\newtheorem{thm}{Theorem}[section]
\newtheorem{lem}[thm]{Lemma}
\newtheorem{cor}[thm]{Corollary}
\newtheorem{df}[thm]{Definition}
\newtheorem{pro}[thm]{Proposition}
\newtheorem{sas}[thm]{Standing Assumption}
\theoremstyle{remark}
\newtheorem{rem}[thm]{Remark}
\title{Actions of $\Z^k$ associated to higher rank graphs}
\author{Alex Kumjian}
\address{Department of Mathematics, University of Nevada, Reno NV
89557, USA}
\email{alex@unr.edu}
\author{David Pask}
\address{Mathematics, SMPS, % and Physical Sciences,
The University of Newcastle, NSW 2308, Australia}
\email{davidp@maths.newcastle.edu.au}
\thanks{Research supported by the Australian Research Council and the
National Science Foundation}
\date{\today}%. Preliminary Draft: Please do not circulate
\begin{document}
\begin{abstract}
An action of $\Z^k$ is associated to a higher rank graph $\Lm$
satisfying a mild assumption.  This generalises the construction of
a topological Markov shift arising from a nonnegative integer matrix.
We show that the stable Ruelle algebra of $\Lambda$ is strongly
Morita equivalent to $C^*(\Lm)$.
% The associated stable Ruelle algebra $R_s$ is
% shown to be strongly Morita equivalent to $C^* (\Lm)$.
Hence, if $\Lambda$ satisfies the aperiodicity condition, the stable
Ruelle algebra is simple, stable and purely infinite.
\end{abstract}
\maketitle
\section{Introduction}\label{intro}
The shift map defines a homeomorphism on the space of two-sided
infinite paths in a finite directed graph, a compact zero dimensional
space when endowed with the natural topology.
% The set of two-sided infinite paths in a finite directed graph may be
% given a natural topology so that it becomes a zero dimensional
% compact space endowed with a distinguished homeomorphism, the shift
% map.
% This dynamical system has provided a key class of examples in
% the field of symbolic dynamics, called Markov shifts or shifts of
% finite type.
Such dynamical systems, called topological Markov shifts or shifts of
finite type, form a key class of examples in symbolic dynamics.
Higher dimensional analogs which exhibit many of the
same dynamical properties include axiom A diffeomorphisms
studied by Smale \cite{sm}. The local hyperbolic nature of the
homeomorphisms in many of the examples has led to the axiomatization
of Smale spaces (see \cite{ru1}). In \cite{ru2,pt1} (see
also \cite{kps} for a short survey and \cite{pt2} for extended notes)
certain $C^*$-algebras were associated to a Smale space making use of the
asymptotic, stable and unstable equivalence relations engendered by
the homeomorphism. The Ruelle algebras, crossed products of the stable
and unstable algebras by the canonical automorphism may be regarded as
higher-dimensional generalisations of Cuntz-Krieger algebras (see
\cite{ck,pt2,pts}). If the graph is irreducible,  the stable Ruelle algebra
associated to the Markov shift is strongly Morita equivalent to the
Cuntz-Krieger algebra associated to the incidence matrix of the graph
(cf.\ \cite[Theorem 3.8]{ck} and \cite[Proposition 3.7]{kps} for
similar results).
%  see \cite[?]{pt2},

% Motivated by the geometrical examples of Robertson and Steger arising
% from group actions on buildings (see \cite{rst1, rst2}), the authors defined
% % a structure $( \Lambda , d)$, called
% a $k$-graph, a higher rank analog of a directed graph (see \cite{kp}).

Following \cite{kp} a $k$-graph is defined to be a higher rank analog
of a directed graph. The definition of a $k$ graph is motivated by the
geometrical examples of Robertson and Steger
arising from group actions on buildings (see \cite{rst1, rst2}).
Given a $k$-graph $\Lambda$, we define a universal $C^*$-algebra,
$C^*( \Lambda )$, the Cuntz-Krieger algebra of $\Lambda$.
Under a mild assumption we form the ``two-sided path space'' of
$\Lambda$, a natural zero dimensional space associated with a
$k$-graph on which there is a $\Z^k$ action by an analog of the shift.
We establish that the key
dynamical properties identified by Ruelle (see \cite{ru1}), when
properly interpreted, hold for this action. Our program then follows the
one set out by Putnam. If $\Lambda$ is irreducible and has finitely many
vertices, then as in \cite{pt1, pt2}, we construct $C^*$-algebras
from the stable and unstable equivalence relations on which there
are natural $\Z^k$ actions. We then form the resulting
crossed products, the %so-called
Ruelle algebras, $R_s$ and
$R_u$. Furthermore, we show that the Ruelle algebra $R_s$ is strongly
Morita equivalent to $C^* ( \Lambda )$. Then, if $\Lambda$
satisfies the aperiodicity condition, the Ruelle algebra
$R_s$ is a Kirchberg algebra, that is, $R_s$ is simple, nuclear and purely
infinite (a similar result holds for $R_u$).
See \cite{pts} for general results on the Ruelle algebras of Smale spaces.
% For related results for Smale spaces see \cite{pts}.

The paper is organised as follows. In section \ref{prelim} we
establish our notation and collect facts for later use. We define a
$k$-graph $( \Lambda , d )$ to be a small category $\Lambda$ equipped with a
degree map $d$ satisfying a certain factorisation property. When
$( \Lambda ,  d)$ satisfies the standing assumption, every vertex of
$\Lambda$ receives and emits a finite but non-zero number of edges of
any given degree, we form $\Lambda^\Omega$ the one-sided infinite path
space of $\Lambda$. Pairs of shift-tail equivalent paths in $\Lambda^\Omega$
% has a topology generated by cylinder sets. For $p \in \N^k$ the shift
% $\sigma^p : \Lambda^\Omega \to \Lambda^\Omega$ is
% a local homeomorphism and
give rise to elements in the path groupoid  $\mathcal{G}_\Lambda$.
The groupoid $C^*$-algebra $C^* ( \mathcal{G}_\Lambda )$ is
naturally isomorphic to $ C^* ( \Lambda )$  (see
\cite[Corollary 3.5]{kp}). There is a canonical
gauge action $\alpha$ of $\T^k$ on $C^* ( \Lambda )$ whose fixed point algebra
$C^* ( \Lambda )^\alpha$ is an AF algebra which coincides with the
$C^*$-algebra of a subgroupoid  $\Gamma_\Lambda$ of $\mathcal{G}_\Lambda$
under this identification. We conclude the section by
stating some facts about principal proper
groupoids.

In section \ref{zkactions} we build a topological dynamical system from a
$k$-graph which is generated by $k$ commuting homeomorphisms. We show
that it satisfies analogs of the two conditions
(SS1) and (SS2) for a Smale space defined in \cite[\S 7.1]{ru1}.
The two-sided path space $\Lambda^\Delta$ of $\Lambda$ has a
zero dimensional topology generated by cylinder sets which is also
given by a metric $\rho$; it is compact if
$\Lambda^0$ is finite. For $n \in \Z^k$ the shift
$\sigma^n : \Lambda^\Delta \to \Lambda^\Delta$  gives
rise to an expansive $\Z^k$-action which is topologically
mixing if $\Lambda$ is primitive.
We show that condition (SS1) is
satisfied, in particular there is a map $(x, y) \mapsto [x, y]$,
defined for $x, y \in \Lambda^\Delta$ with $\rho ( x , y ) < 1$ taking
values in $\Lambda^\Delta$,
% $$
% [ \cdot , \cdot ] : \{ ( x,y ) \in \Lambda^\Delta
% \times \Lambda^\Delta : \rho ( x , y ) < 1 \} \to \Lambda^\Delta ,
% $$
which endows the space $\Lambda^\Delta$ with a local product
structure.  For $x \in \Lambda^\Delta$ there are subsets $E_x$ and
$F_x$ of $\Lambda^\Delta$ such that $E_x \times F_x$ is homeomorphic
to a neighbourhood of $x$ (under this bracket map). Moreover, if
$e = ( 1 , \ldots, 1 ) \in \Z^k$ then the shift
$\sigma^e$ contracts the distance between points in $E_x$ and expands
them on $F_x$.  This is our analog of condition (SS2) for a single
homeomorphism. As in \cite{pt1} we define the stable and unstable
relations which may be characterised in terms of tail equivalences on
$\Lambda^\Delta$, since the topology of $\Lambda^\Delta$ is generated by
cylinder sets. The stable and unstable relations
give rise to the stable and unstable groupoids, $G_s$ and $G_u$. Since
the unstable relation for $\Lambda^\Delta$ is exactly the stable relation for
the opposite $k$-graph $\Lambda^{\text{op}}$ (the $k$-graph formed by
reversing all the arrows of $\Lambda$), we focus our attention on the
stable case. Finally, we examine the internal structure of the stable
groupoid $G_s$; it is the inductive limit of a sequence of mutually
isomorphic principal proper groupoids $G_{s,m}$, for $m \in \Z^k$.

In section \ref{ruelle} we associate certain $C^*$-algebras to
an irreducible $k$-graph $\Lambda$ with $\Lambda^0$ finite.
First we state a suitable version of the Perron-Frobenius theorem,
which gives rise to a shift-invariant measure $\mu$ on $\Lambda^\Delta$.
% (see Proposition \ref{mudef}).
The measure $\mu$ decomposes in a manner which respects the local
product structure; this in turn gives rise to Haar systems for $G_s$
and $G_u$.  % (see Proposition \ref{haarhaar}).
The stable and unstable $C^*$-algebras may then be defined:
$S = C^* ( G_s )$ and $U = C^* ( G_u )$.
The $\Z^k$-action on $\Lambda^\Delta$ induces actions $\beta_s$ on $S$
and $\beta_u$ on $U$ which  % (by Proposition \ref{mutrace})
scale the canonical densely-defined traces. % defined using $\mu$.
The Ruelle algebras are defined to be the corresponding
crossed products, $R_s = S \times_{\beta_s} \Z^k$ and
$R_u = U \times_{\beta_u} \Z^k$.

In the last section we prove our main results:
Suppose that $\Lambda$ is an irreducible $k$-graph which has
finitely many vertices. Then
\begin{itemize}
\item[(i)] $S$ is strongly Morita equivalent to
$C^* ( \Lambda )^\alpha$ (see Theorem \ref{zisequiv}),
\item[(ii)] $R_s$ is strongly Morita equivalent to
$C^* ( \Lambda  )$ (see Theorem \ref{zequiv2}).
\end{itemize}

Similar assertions hold for $U$ and $R_u$ when $\Lambda$ is replaced by
$\Lambda^{\text{op}}$.
% We use the notion of equivalence of
% groupoids (in the sense of \cite{mrw}) to establish these results.
% From these results flow many important consequences:
We establish our main results using the notion of equivalence of
groupoids (in the sense of \cite{mrw}). From established properties of
$C^* ( \Lambda )$ flow many important consequences:
The stable
algebra $S$ is an AF algebra and if $\Lambda$ is primitive then $S$ is
simple. The Ruelle algebra $R_s$ is nuclear and in the bootstrap class
$\mathcal{N}$ for which the UCT holds. Further, if $\Lambda$ satisfies the
aperiodicity condition then $R_s$ is simple, stable and purely infinite. The
Kirchberg-Phillips theorem therefore applies, so the
isomorphism class of $R_s$ is completely determined by its $K$-theory
(see \cite{ki,ph}). \\[-3mm]

{\em Acknowledgments}:
The first author wishes to thank Valentin Deaconu for some useful
discussions at an early stage of this work.
He also wishes to thank the second author and his colleagues
at the School of Mathematics and Physical Sciences of the University
of Newcastle for their warm hospitality. % and support

%, where most of this work was done,
\section{Preliminaries}\label{prelim}

In this section we first give a little background, then we establish our
notation and conventions about a $k$-graph $\Lambda$ and its path
groupoid $\mathcal{G}_\Lambda$ which are taken from \cite{kp}. We
define the $C^*$-algebra of a $k$-graph, $C^* ( \Lambda )$, which may
be realised as $C^* ( \mathcal{G}_\Lambda )$. Finally, we
state some results concerning principal proper groupoids and their
Haar systems which are taken from \cite{rn2, mw1, kmrw}.

We shall use $\N$ to denote the set of natural numbers
$\{ 0 ,  1 , 2 , \ldots  \}$; $\Z$, $\R$, $\T$ denote the sets of integers,
real numbers and complex numbers with unit modulus, respectively. For
$k >0$ we endow $\N^k$ and $\Z^k$ with the coordinatewise ordering.

A category is said to be small if its morphisms form a set; the objects
are often identified with a subset of morphisms ($x \mapsto 1_x$).
A {\em groupoid} is a small category $\Gamma$ in which every morphism is
invertible. For $\gamma \in G$ we have $r ( \gamma ) = \gamma \gamma^{-1}$ and
$s ( \gamma ) = \gamma^{-1} \gamma$, then $r,s :  \Gamma \to \Gamma^0$ where
$\Gamma^0$ is the {\em unit space} (or space of objects) of
$\Gamma$. If the groupoid $\Gamma$ is furnished with a topology for
which the groupoid operations are continuous then $\Gamma$ is called a {\em
topological groupoid}.  We shall assume that our groupoids
are equipped with a locally compact, Hausdorff, second
countable topology. If the groupoid $\Gamma$ has a {\em left Haar system}
$\mu = \{ \mu^x : x \in \Gamma^0 \}$, an equivariant system of
measures on the fibres $r^{-1} (x)$, then we may form the full and
reduced $C^*$-algebras, $C^* ( \Gamma )$ and $C^*_r ( \Gamma )$. Since
we shall only be dealing with left Haar systems we shall henceforth
omit the qualifier left.  If $\Gamma$ is amenable then its
full and reduced $C^*$-algebras coincide. The groupoid $\Gamma$ is
called {\em $r$-discrete} if $r$ is a local homeomorphism; in this case the
counting measures form a Haar system.  For more definitions and
properties of groupoids and their $C^*$-algebras, consult
\cite{rn1,m}. For the most part we have followed the
conventions of \cite{rn1}, with the exception that $s$ replaces $d$
for the source map. A good reference for amenable groupoids may be found
in \cite{ar}. We shall frequently invoke the notion of equivalence of groupoids
(\cite[Definition 2.1]{mrw}) which (in the presence of Haar systems)
gives rise to the strong Morita equivalence of their $C^*$-algebras
(\cite[Theorem 2.8]{mrw}).  A good reference for $C^*$-algebras and
their crossed products is \cite{pd}.

Let $k$ be a positive integer.  Recall the notion of $k$-graph (see
\cite{kp}).
\begin{df}
A $k$-graph is a pair $(\Lm, d)$, where $\Lm$ is a countable small
category and $d : \Lm \to \N^k$ is a morphism, called the degree map,
such that the factorisation property holds: for every
$n_1, n_2 \in \N^k$ and $\lm \in \Lm$ with $d(\lm) = n_1 + n_2$, there exist
unique elements $\nu_1 , \nu_2 \in \Lm$ with
$$
\lm = \nu_1 \nu_2, \quad n_1 = d (\nu_1 ), \quad n_2 = d( \nu_2 ).
$$
For $n \in \N^k$ write $\Lm^n = \{ \lm \in \Lm : d(\lm) = n \}$.  It
will be convenient to identify $\Lm^0$ with the objects of $\Lm$.
Let $r, s: \Lm \to \Lm^0$ denote the range and source maps.
\end{df}

Let $E = (E^0, E^1)$ be a (countable) directed graph.  Then the set of
finite paths $E^*$ together with the length map defines a $1$-graph
(the roles of $r$ and $s$ must be switched).

If $\Lambda$ is a $k$-graph then the opposite category
$\Lambda^{\text{op}}$ can also be made into a $k$-graph by setting
$d ( \lambda^{\text{op}} ) = d ( \lambda )$.
% denote the opposite $k$-graph by $\tilde{\Lambda}$.

The $k$-graph which gives us the prototype for a (one-sided) infinite
path is
$$
\Omega = \Omega_k = \{ (m,n) : m,n \in \N^k : m \le n \} .
$$
The structure maps are given by
\begin{equation} \label{structuremaps}
r(m,n)=m ,\quad s (m,n)=n, \quad ( \ell,n) = ( \ell,m) (m,n) , \quad
d ( m , n ) = n - m
\end{equation}
%$\Omega$ is a small category
where the object space is identified with $\N^k$ % .
% Define   $d : \Omega \rightarrow \N^k$ by $d (m,n)= n-m$
(see \cite[Example 1.7ii]{kp}).
For other examples of $k$-graphs
consult \cite{kp}. % and \cite{rsy}.

\begin{df}
A $k$-graph $\Lambda$ is said to be irreducible (or strongly connected)
if for every $u, v \in \Lm^0$, there is $\lm \in \Lm$ with
$d(\lambda) \ne 0$ such that $u = r(\lm)$ and $v = s(\lm)$.
We say that $\Lambda$ is primitive if there is a nonzero $p \in \N^k$  so that
for every $u , v \in \Lambda^0$ there is $\lambda \in \Lambda^p$ with
$r ( \lambda ) = u$ and $s ( \lambda )=v$.
\end{df}

Suppose that $\Lambda$ is primitive; then there is an $N$ such that for
all $p \ge N$ and every $u , v \in \Lambda^0$ there is
$\lambda \in \Lambda^p$ with $r ( \lambda ) = u$ and $s ( \lambda
)=v$. Moreover, under the following standing assumption \ref{star}, $\Lambda^0$
must be finite.

To ensure that the analog of the two-sided infinite path space (to be
discussed in the next section) is nonempty and locally compact we
shall need the following standing hypothesis.

\begin{sas} \label{star}
For each $p \in \N^k$ the restrictions of $r$ and $s$ to $\Lm^p$ are
surjective and finite to one.
\end{sas}

The standing hypothesis used here is equivalent to the requirement
that both $\Lambda$ and $\Lambda^{\text{op}}$ satisfy the condition of
\cite[\S 1]{kp}. Recall from \cite{kp} the definition of the universal
$C^*$-algebra of a $k$-graph.

\begin{df} \label{cstarlambdadef}
Let $\Lambda$ be a  $k$-graph. Then $C^* ( \Lambda )$ is
defined to be the universal $C^*$-algebra generated by a family
$\{ s_\lambda : \lambda \in \Lambda \}$ of partial isometries
satisfying:
\begin{itemize}
\item[(i)] $\{ s_v : v \in \Lambda^0 \}$ is a family of mutually
orthogonal projections,
\item[(ii)] $s_{\lambda \mu} = s_\lambda s_\mu$ for all
$\lambda , \mu \in \Lambda$ such that $s ( \lambda ) = r ( \mu )$,
\item[(iii)] $s_\lambda^* s_\lambda = s_{s( \lambda )}$ for all
$\lambda \in \Lambda$,
\item[(iv)] for all $v \in \Lambda^0$ and $n \in {\bf N}^k$ we  have
$\displaystyle s_v = \sum_{\lambda \in
\Lambda^n ( v )} s_\lambda s_\lambda^*$.
\end{itemize}
\end{df}

Let $\Lambda$ be a $k$-graph and set
\begin{equation} \label{lambdaomega}
\Lambda^\Omega = \{ x : \Omega \rightarrow \Lambda : x \text{ is a
$k$-graph homomorphism } \} ;
\end{equation}
note that a $k$-graph morphism must preserve degree, so that for
$x \in \Lambda^\Omega$, we have
$d(x(m,n)) = n-m$. In \cite{kp} the set $\Lambda^\Omega$ was
denoted $\Lambda^\infty$.
By \ref{star}, $\Lm^\Omega \ne \nul$. For each $\lambda \in \Lambda$ we
put
\begin{equation} \label{topdef}
Z ( \lambda ) = \{ x \in \Lambda^\Omega : x  ( 0 ,  d(\lambda ) ) =
\lambda \}
\end{equation}
then again by \ref{star}, $Z( \lambda )\neq \nul$. The collection
of all such cylinder sets forms a basis for a topology on $\Lambda^\Omega$
under which each such subset is compact. For $p \in \N^k$ define a map
$\sigma^p : \Lambda^\Omega \rightarrow \Lambda^\Omega$ by
\begin{equation} \label{sigmadef}
( \sigma^p x ) (m,n) = x ( m+p , n + p )
\end{equation}
note that $\sigma^p$ is a local homeomorphism. Now we
form the {\em path groupoid} (for more
details see \cite{kp})
$$
\mathcal{G}_\Lambda = \{ (x , n , y ) : x , y \in \Lambda^\Omega , n
\in \Z^k ,
\sigma^\ell x = \sigma^m y , n = \ell -m \text{ for some } \ell , m
\in \N^k \} ,
$$
with structure maps
$$
r ( x,n,y) = x , \quad s ( x,n,y) = y , \quad \text{and} \quad
(x,m,y)(y,n,z) = (x,m+n,z) ,
$$
where we have identified $\Lambda^\Omega$ with the unit space by
$x \mapsto (x,0,x)$.
% The topology on
% $\mathcal{G}_\Lambda$ is generated by compact open sets
% \begin{equation} \label{ggroupoidtop}
% Z ( \lambda , \mu ) = \{ (x, d ( \lambda ) - d ( \mu ) ,y) :
% x ( 0 ,  d (\lambda ) ) = \lambda , y ( 0 , d( \mu ) ) = \mu ) \}
% \end{equation}
% for
% $\lambda , \mu \in \Lambda$ with $s ( \lambda ) = s ( \mu)$.
% One of
% the difficulties of working with $\mathcal{G}_\Lambda$ is that is does
% not obviously arise as a semidirect product groupoid for a $\Z^k$
% action. Following the methods of \cite{pt1, pt2} one way around this is to
% consider a two-sided version of the above:

By  \cite[Corollary 3.5(i)]{kp}
$C^* ( \mathcal{G}_\Lambda ) = C^* ( \Lambda )$.
There is a canonical gauge action
$\alpha : \T^k \to \Aut (C^* ( \Lambda))$ which is realised on the
dense subalgebra $C_c (\mathcal{G}_\Lambda )$ by
$$
\alpha_t ( f ) ( x,n,y) = t^{n} f(x,n,y)
$$
where $t^{n} = \prod_i t_i^{n_i}$.
The fixed point algebra for this action $C^* ( \Lambda )^\alpha$ is the
closure of the subalgebra of $C_c ( \mathcal{G}_\Lambda )$ consisting
of functions which vanish at points of the form $(x,n,y)$ with
$n \neq 0$. Hence $C^* ( \Lambda )^\alpha$ is isomorphic
to $C^* ( \Gamma_\Lambda )$ where $\Gamma_\Lambda$ is the open
subgroupoid of $\mathcal{G}_\Lambda$ given by
$$
\Gamma_\Lambda = \{ (x ,0, y ) : x , y \in \Lambda^\Omega ,
\sigma^m x = \sigma^m y \mbox{ for some } m \in \N^k \} .
$$
Recall that
$$
C^* ( \Lambda )^\alpha = \mathcal{F}_\Lambda = \lim_{m \rightarrow \infty}
\mathcal{F}_{m}
$$
where
$$
\mathcal{F}_m \cong \bigoplus_{v \in \Lambda^0} \mathcal{K} \left(
 \ell^2 ( \{ \lambda \in \Lambda^m : s ( \lambda ) = v \} ) \right) ;
$$
hence, $C^* ( \Gamma_\Lambda ) = C^* ( \Lambda )^\alpha$
is an AF algebra (see \cite[Lemmas 3.2, 3.3]{kp}).

Let $\Lambda_i$ be a $k_i$ graph for $i=1,2$, then
$\Lambda_1 \times \Lambda_2$ is a $(k_1 + k_2 )$-graph in a natural
way (see \cite[Proposition 1.8]{kp}). By \cite[Corollary 3.5iv]{kp} we have
$$
C^* ( \Lambda_1 \times \Lambda_2 ) \cong
C^* ( \Lambda_1 ) \otimes C^*  ( \Lambda_2 ) .
$$
If $k_1 = k_2 =k$ then we may form a $k$-graph
$$
\Lambda_1 \diamond \Lambda_2 =
\{ ( \lambda_1 , \lambda_2 ) \in \Lambda_1 \times \Lambda_2 :
d ( \lambda_1 ) = d ( \lambda_2 ) \}
$$
with $d( \lambda_1 , \lambda_2 ) = d( \lambda_1 )$ and the other
structure maps inherited from $\Lambda_1 \times \Lambda_2$. Note that
$\Lambda_1 \diamond \Lambda_2 = f^*(\Lambda_1 \times \Lambda_2 )$
where $f : \N^k \to  \N^k \times \N^k$ is given by
$f(m) = (m, m)$ (cf.\ \cite[Example 1.10iii]{kp}).

Let $G$ be a compact abelian group and for $i=1,2$ let
$\alpha^i : G \to \Aut(A_i)$ be a strongly continuous action of $G$ on the
$C^*$-algebra $A_i$. Let $A_1 \otimes_G A_2$ denote the
fixed-point algebra $( A_1 \otimes A_2 )^\eta$ where
$\eta : G \to \Aut (A_1  \otimes A_2 )$ is given by
$\eta_g(a \otimes b) = \alpha^1_g(a) \otimes \alpha^2_{g^{-1}}(b)$.
This is the natural notion of tensor product in the category of
$C^*$-algebras with a given $G$-action see \cite[\S 2]{opt}.
Now with $\Lambda_i$ as above and taking $\alpha^i$ to be the gauge
action on $C^* ( \Lambda_i )$ for $i=1,2$ we have
\begin{equation} \label{diamond}
C^* ( \Lambda_1 \diamond \Lambda_2 ) \cong
C^* ( \Lambda_1 ) \otimes_{\T^k} C^* ( \Lambda_2 ).
\end{equation}
This follows by an argument similar to the proof of
\cite[Proposition 2.7]{ku2}.
 %Let $\alpha^i$ denote the gauge action on $C^* ( \Lambda_i )$ for $i=1,2$.

In the remainder of this section we state some standard facts
concerning principal proper groupoids in a convenient form.
Recall that a groupoid $G$ is said to be {\em principal},
if it is isomorphic to an equivalence relation, that is, if
$r \times s : G \to G^0 \times G^0$ is an embedding. If, in addition,
the image is a closed subset of $G^0 \times G^0$, it is said
to be  {\em proper} (see \cite{mw1}).
\begin{lem} \label{ppg}

Let $\pi : X \rightarrow Y$ be a continuous open surjection  between two
locally compact Hausdorff spaces. Then
$$
X \star_\pi X = \{ (x,y) \in X \times X : \pi (x) = \pi (y) \}
$$
is a principal proper groupoid, with structure maps $r(x,y)=x$,
$s(x,y)=y$ and $(x,y)(y,z)=(x,z)$.
Moreover, $X$ is an $( X \star_\pi X , Y )$-equivalence.
% Moreover, every principal proper groupoid is of this form.
\end{lem}
\begin{proof}
Evidently $X \star_\pi X$ is a principal groupoid; since
$X \star_\pi X$ is a closed subset of $X \times X$ it is proper.
By \cite[Example 2.5]{mrw} $X$ is a $( X \star_\pi X , Y )$-equivalence.
% Let $G$ be a principal proper groupoid, then $g \mapsto ( s(g),r(g))$ is a
% homeomorphism onto a closed subset of $G^0 \times G^0$. It follows that
% $G^0 /  G$ is Hausdorff and the quotient map
% $\pi : G^0 \to G^0 / G$ is continuous and open.
\end{proof}

The following definition is taken from \cite[\S 1]{rn2} (see also
\cite[Definition 5.42]{m}).

\begin{df}
With $\pi$ be as above, a $\pi$-system consists of a family
$$
\mu = \{ \mu^y : y \in Y \}
$$
of positive Radon measures on $X$ such that the support of $\mu^y$ is
contained in $\pi^{-1} (y)$ for each $y \in Y$ and the function
$$
\mu ( f ) (y) = \int f (x) d \mu^y (x)
$$
lies in $C_c (Y)$ for each $f \in C_c (X)$. If the support of each
$\mu^y$ is all of $\pi^{-1} (y)$ for all $y \in Y$, then the
$\pi$-system is said to be full.
\end{df}

Note that a Haar system on a groupoid is an equivariant $r$-system.
A full $\pi$-system gives rise to a Haar system for $X \star_\pi X$.

\begin{pro} \label{pisys}
Let $\pi$  be as above and $\mu = \{ \mu^y : y \in Y \}$
be a full $\pi$-system. Then for $x \in X$
$$
\tilde{\mu}^x = \delta_x \times \mu^{\pi (x)}
$$
defines a Haar system $\tilde{\mu} = \{ \tilde{\mu}^x : x \in X \}$
for $X \star_\pi X$. Moreover, $C^* ( X \star_\pi X )$ is strongly
Morita equivalent to $C_0 ( Y )$. There is a densely defined
$C_0 (Y)$-valued trace on $C^* (X \star_\pi X )$ given by
\begin{equation} \label{trace}
\tau_\mu (f) (y) = \int_X f (x,x) d \mu^y (x)
\end{equation}
for $f \in C_c ( X \star_\pi X )$.
\end{pro}

\begin{proof}
%% The existence of the Haar system
The first assertion follows from \cite[Proposition
5.2]{kmrw} (see also \cite[Theorem 5.51]{m}). The Morita equivalence
now follows from \cite[Theorem 2.8]{mrw} and Lemma \ref{ppg}
(see also \cite[Proposition 2.2]{mw1}). A routine computation
shows that $\tau_\mu (fg) = \tau_\mu (gf)$ for
$f,g \in C_c ( X  \star_\pi X )$.
\end{proof}

\section{$\Z^k$ actions}\label{zkactions}

In this section we adapt the methods of \cite{pt1} for the $\Z$-action
associated to shift of finite type to analyze an analogous $\Z^k$
action on a topological space associated to a $k$-graph. Many of the
constructions of \cite{ru1, pt1} can be generalised to this
setting. Following \cite{pt1, pt2}, we provide a
description of the stable, unstable and asymptotic relations for our
$\Z^k$ action and the topology of the associated groupoids.

There is a natural $\Z^k$ action on the analog of the two-sided path
space of a $k$-graph satisfying the standing hypothesis \ref{star}.
First, we form a $k$-graph
which gives us the prototype
of a two-sided infinite path:
Set
$$
\Dl = \Dl_k = \{ (m, n) : m, n \in \Z^k, m \le n \} ;
$$
with % $d :\Dl \to \N^k$ given by $d(m, n) = n - m$ and other
% We identify $\Z^k$
% with $\Dl^0$ in the obvious way ($n \mapsto (n, n)$).  The structure
% maps of $\Dl$ are given as follows:
% $$
% r(m, n) = m, \quad s(m, n) = n, \quad (l, n) = (l, m)(m, n).
% $$
structure maps given as in (\ref{structuremaps}), it is
straightforward to check that $( \Delta , d )$ is a $k$-graph.
% $\Delta$ is a small category and that the factorization property is
% satisfied. Hence
% We now construct a topological space canonically associated with a
% $k$-graph $\Lambda$ in the same way the space of (two-sided) infinite
% sequences is associated to a matrix with nonnegative integer matrices
% (it is usually assumed that there are no rows or columns consisting
% exclusively of zeroes).

Next we use $\Delta$ to form the two-sided infinite path space (cf.\
(\ref{lambdaomega})). Set
$$
\Lm^\Dl = \{ x: \Dl \to \Lm : x \text{ is a $k$-graph morphism} \} ;
$$
% note that a $k$-graph morphism is degree preserving (i.e.\
% $d(x(\dl)) = d(\dl)$).
then by \ref{star}, $\Lm^\Dl \ne \nul$.
We endow $\Lm^\Dl$ with a topology as follows (cf.\ (\ref{topdef})):
for each
$n \in \Z^k$ and $\lm \in \Lm$ set % with $r(\lm) = n$
\begin{equation*} \label{tstopdef}
Z(\lm, n) =  \{ x \in \Lm^\Dl : x(n, n + d(\lm)) = \lm \}.
\end{equation*}
Again by \ref{star}, $Z(\lm, n)  \ne \nul$.  The collection of all
such cylinder sets forms a basis for a topology on $\Lm^\Dl$ for which
each such subset is compact.  It follows that $\Lambda^\Delta$ is a
zero dimensional space and if $\Lm^0$ is finite,
then $\Lm^\Dl$ is itself compact (since
$\Lm^\Dl = \cup_{v \in \Lm^0} Z(v, 0)$).
Now for each $n \in \Z^k$ we define a map
$\sg^n : \Lm^\Dl \to \Lm^\Dl$  by
$$
\sg^n (x)(\ell , m ) = x(\ell + n, m + n) .
$$
% The reader may wish to check that $\sg^p$ is a homeomorphism.
% Moreover
Note that $\sigma^n$ is a homeomorphism for every $n \in \Z^k$,
$\sg^{n + m} = \sg^n\sg^m$  for $n, m \in \Z^k$ and
$\sigma^0$ is the identity map.

We define a metric on $\Lm^\Dl$ as follows.
We set $e = ( 1 , \ldots , 1 ) \in \Z^k$ and for $j \in \N$, let
$\theta_j \in \Dl$ denote the element $(-je, je)$;
note that $\theta_0 = 0$.  Given $x, y \in \Lm^\Dl$, set
$$
h(x, y) =
\begin{cases}
0 & x(0) \ne y(0) \\
1 + \sup \{j : x(\theta_j) = y(\theta_j) \} & \text{ otherwise.}
\end{cases}
$$
Fix $0 < r < 1$; we may define a metric $\rho$ on $\Lm^\Dl$ by the formula
$\rho(x, y) = r^{h(x, y)}$ for $x, y \in \Lm^\Dl$ (note that
$\rho (x,x) = r^\infty = 0$).
The topology induced by this metric is the same as the one above.

\begin{pro} \label{dynfacts}
The $\Z^k$-action
$n \mapsto \sigma^n$ on $\Lambda^\Delta$ is  expansive in the sense
that there is an $\varepsilon > 0$ such that for all
$x, y \in \Lambda^\Delta$ if
$\rho ( \sigma^n (x) , \sigma^n (y ) ) < \varepsilon$
for all $n$ then $x=y$. Moreover, if $\Lambda$ is primitive then $\sigma$ is
topologically mixing in the sense that for any two nonempty open
sets $U$ and $V$ in $\Lambda^{\Delta}$ there is a $Q \in \Z^k$ so
that $U \cap \sigma^q(V) \ne \nul$ for all $q \ge Q$.
\end{pro}

\begin{proof}
To show that the action is expansive, observe that $\varepsilon = r$
will suffice (if $ x( n-e,n+e ) = y (n-e,n+e)$ for all $n \in \Z^k$,
then $x=y$). If $\Lambda$ is primitive there is an $M \in \N^k$ such that for
all $m \ge M$ and every $u , v \in \Lambda^0$ there is
$\lambda \in \Lambda^m$ with $r ( \lambda ) = u$ and
$s ( \lambda )=v$. To show that for any two nonempty open
sets $U$ and $V$ in $\Lambda^{\Delta}$ there is a $Q \in \Z^k$ so
that $U \cap \sigma^q(V) \ne \nul$ for all $q \ge Q$, it suffices
to demonstrate this for cylinder sets. So let
$U =  Z ( \lambda , \ell )$ and $V =Z ( \nu , n )$.
Set $Q = M  + d ( \nu ) + n - \ell$; then given $q \ge Q$,
there is $\lambda' \in \Lambda$ with $d ( \lambda' ) = M + q - Q$ such
that $r ( \lambda' ) = s ( \nu )$ and
$s ( \lambda' ) = r   ( \lambda )$. Observe that
$$
Z (\nu \lambda'\lambda , n-q ) \subset Z ( \lambda , \ell ) \cap \sigma^q
( Z ( \nu , n ) ).
$$
Therefore,  $U \cap \sigma^q (V) \ne \nul$ for all $q \ge Q$, as required.
\end{proof}

\begin{rem}\label{restriction}
Consider the $1$-graph obtained by restricting consideration
to powers of $\Lambda^e$; note that this $1$-graph may be regarded
as the $k$-graph $f^*(\Lambda)$ where $f: \N \ra \N^k$ is given by
$f(j) = je$ (see \cite[Definition 1.9]{kp}).  By arguing as in
\cite[Proposition 2.9]{kp} it follows that the restriction map
$\Lambda^{\Delta} \ra f^*(\Lambda)^{\Delta}$ is a homeomorphism.
Under this identification the generator of the action of $\Z$ on
$f^*(\Lambda)^{\Delta}$ is identified with $\sg^e$. 
Many attributes of this restricted dynamical system are reflected in
the action of $\Z^k$, as we shall see below.
%% This restricted action contains a great deal of inherent dynamical
%% information, as we shall see when we consider the contracting and
%% expanding directions below.
\end{rem}

The space $\Lambda^\Delta$ decomposes locally into contracting and
expanding directions for the shift. For $x \in \Lambda^\Delta$ set
\begin{align*}
E_x &= \{ y \in \Lambda^\Delta :
 x ( m,n ) = y ( m , n ) , \mbox{ for all } 0 \le m \le n \} \\
F_x &= \{ y \in \Lambda^\Delta : x ( m , n ) = y ( m , n ) ,
\mbox{ for all } m \le n \le 0 \} .
\end{align*}

Observe that for $j \in \N$ we have (see \cite[\S 7.1]{ru1}, also
\cite{pt1})
\begin{align*}
\rho ( \sigma^{je} (y) , \sigma^{je} (z) ) &\le r^j \rho ( y,z )
\mbox{ for } y,z \in E_x \\
\rho ( \sigma^{-je} (y) , \sigma^{-je} (z) ) &\le r^j \rho ( y,z )
\mbox{ for } y,z \in F_x .
\end{align*}
For $p \ge 0$ a simple calculation shows
that $\sigma^p E_x \subseteq E_{\sigma^p x}$ and
$\sigma^{-p} F_x \subseteq F_{\sigma^{-p} x}$.

\begin{pro} \label{bracket} \textnormal{(cf.\ \cite{ru1, pt1})}
There exists a unique map
$$
[ \cdot , \cdot ] : \{ (x,y) \in \Lambda^\Delta \times \Lambda^\Delta
: \rho (x,y) < 1 \} \to \Lambda^\Delta
$$
satisfying
\begin{equation} \label{brackdef}
\begin{aligned}[]
[x,y](m,n) &= x ( m , n ) \quad \text{if } m \le n \le 0  \\
[x,y](m,n) &= y ( m , n ) \quad \text{if } 0 \le m \le n  .
\end{aligned}
\end{equation}
Moreover, $[ \cdot , \cdot ]$ is continuous,
$F_x \cap E_y = \{ [x,y] \}$ if $\rho (x,y)<1$ and the following hold
\begin{equation} \label{brackprops}
[x,x] = x  , \quad [[x,y],z]=[x,z] , \quad  [x,[y,z]] = [x,z] , \quad
[ \sigma^n x , \sigma^n y ] = \sigma^n [x,y] ,
\end{equation}
wherever both sides of each equation are defined. Furthermore, for
$x \in \Lambda^\Delta$ the
restriction of $[ \cdot , \cdot ]$ to $E_x \times F_x$ induces a
homeomorphism $E_x \times F_x \cong Z(x(0),0)$.
\end{pro}

\begin{proof}
By the  factorisation property, a consistent family of elements
$$
\{ x (-p,p ) \in \Lambda^{2p} : p \ge 0 \} \text{ with }
x ( -q , q ) = \lambda x (-p,p) \mu
$$

for some $\lambda , \mu$ when
$q \ge p$, will determine a unique element $x \in \Lambda^\Delta$
(cf.\ \cite[Remarks 2.2]{kp}).
For $p \ge 0$ and $x,y$ with $\rho (x,y)<1$ (so that $x(0)=y(0)$) set
$$
[x,y](-p,p) = x(-p,0) y (0,p) .
$$
It is straightforward to check that this results in the unique  map
satisfying (\ref{brackdef}); moreover it is continuous.
If $z \in F_x \cap E_y$, then $z (m,n) = x(m,n)$ for $m \le n \le 0$ and
$z(m,n) = y(m,n)$ for $0 \le m \le n$; hence, $z = [x,y]$ by
(\ref{brackdef}).

The properties (\ref{brackprops}) are straightforward to verify.
For the last assertion, it is clear that the restriction of
$[ \cdot , \cdot ]$ to $E_x \times F_x$ is one-to-one. To see that the
image is $Z(x(0),0)$, let $z \in Z ( x(0), 0)$; then $[z,x] \in E_x$,
$[x,z] \in F_x$ and $z = [[z,x],[x,z]]$. The restriction is clearly
continuous  as is its inverse $z \mapsto ([z,x],[x,z])$.
\end{proof}

Note that if $\rho ( x, y ) < 1$, then $y \in E_x$ if and only if
$[x,y]=x$ and similarly $y \in F_x$ if and only if $[y,x]=x$.

As in \cite{pt1}  we define the stable and unstable equivalence
relations on $\Lm^\Dl$ as follows.  Given $x, y \in \Lm^\Dl$ define
\begin{align*}
x \sim_s y & \quad \text{ if }
\lim_{j \ra \infty} \rho(\sg^{je} (x), \sg^{je} (y)) = 0 \\
x \sim_u y & \quad \text{ if }
\lim_{j \ra -\infty} \rho(\sg^{je}(x), \sg^{je}(y)) = 0.
\end{align*}

Note that $x \sim_s y$ if and only if
there is $m \in \Z^k$ such that for all $n \in \Z^k$ with $m \le n$ we
have $x(m, n) = y(m, n)$. Similarly $x \sim_u y$ if and only if
there is $n \in \Z^k$ such that for all $m \in \Z^k$ with $m \le n$ we
have $x(m, n) = y(m, n)$.

These equivalence relations give rise to two locally compact
groupoids: the stable groupoid,
$$
G_s = G_s ( \Lambda ) =
\{ (x, y) \in \Lm^\Dl \times \Lm^\Dl : x \sim_s y  \}
$$
and the unstable groupoid,
$$
G_u = G_u ( \Lambda ) =
\{ (x, y) \in \Lm^\Dl \times \Lm^\Dl : x \sim_u y  \};
$$
the unit space of each is identified with $\Lm^\Dl$ and the structure
maps are the natural ones.
The topology on $G_s$ is given as follows: For $m \in \Z^k$ set
$$
G_{s, m} = \{ (x, y) \in \Lm^\Dl \times \Lm^\Dl :
               x(m, n) = y(m, n) \text{ for all }n \ge m \}.
$$
Note that $G_{s,m}$ is a subgroupoid of $G_s$.
We endow $G_{s, m}$ with the relative topology
%(note that $G_{s,m}$ is closed)
and $G_s = \cup_m G_{s,m}$ with
the inductive limit topology.  The topology on $G_u$ is defined
similarly. None of these groupoids are $r$-discrete in general.

There is an natural inclusion map $\Omega \hookrightarrow \Delta$ which
gives rise to a surjective map
$\pi : \Lambda^\Delta \rightarrow \Lambda^\Omega$
given by restriction: $\pi ( x ) (m,n) =
x(m,n)$ for $x \in \Lambda^\Delta$ and $(m,n) \in \Omega$.
It is straightforward
to verify that $\pi$ is continuous and open. Observe that for $x \in
\Lambda^\Delta$ and $p \in \N^k$ we have
\begin{equation} \label{stationary}
\pi \circ \sigma^p (x) = \sigma^p \circ \pi (x) .
\end{equation}
We now collect some facts
about the topology of $G_s$ for future use:

\begin{pro} \label{gsprops}
Let $\Lambda$ be a $k$-graph, and $G_s$ the groupoid defined
above. Then for all $m \in \Z^k$, $G_{s,m}$ is a closed subset of
$\Lambda^\Delta \times \Lambda^\Delta$; indeed
$$
G_{s,m} = \Lambda^\Delta \star_{\pi \circ \sigma^m} \Lambda^\Delta = \{ (
x,y) \in \Lambda^\Delta \times \Lambda^\Delta : \pi ( \sigma^m x) =
\pi ( \sigma^m y) \}
$$
and hence $G_{s,m}$ is a
principal proper groupoid.
For all $m, n \in \Z^k$ we have that
$G_{s,m+n} = ( \sigma^{-m} \times \sigma^{-m} ) G_{s,n}$;
in particular, the $G_{s,m}$ are all  isomorphic to
$\Lambda^\Delta \star_\pi\Lambda^\Delta$. Moreover, for $m \le n$,
$G_{s,m}$ is an open subset of $G_{s,n}$.
% For $(x,y) \in \Lambda^\Delta \times \Lambda^\Delta$ we have
% $(x,y) \in G_s$ if and only if
% $( \pi (x) , 0 , \pi (y ) ) \in \Gamma_\Lambda$.
\end{pro}

\begin{proof}
For the first part, observe that
$\pi ( \sigma^m x) = \pi ( \sigma^m y)$ if and only if
$x (m,n) = y(m,n)$ for all $n \ge m$; so
$G_{s,m} = \Lambda^\Delta \star_{\pi\circ\sigma^m} \Lambda^\Delta$
which is a principal proper groupoid by Lemma \ref{ppg}. For the
second assertion, note that
$x(n+m,\ell)=y(n+m,\ell)$ for all $\ell \ge n+m$ if and only if
$\sigma^{m} x (n,\ell') = \sigma^{m} y(n,\ell')$ for all
$\ell' \ge n$ and that
$( \sigma^{-m} \times \sigma^{-m} )$ is a homeomorphism of
$\Lambda^\Delta \times \Lambda^\Delta$. To show that $G_{s,m}$ is an
open subset of $G_{s,n}$ for $m \le n$, it suffices
to consider the case when $m=0$.  Suppose that $(x,y) \in G_{s,0}$, then
we have $x,y \in Z ( \lambda , 0 )$ where
$\lambda = x ( 0 , n ) = y ( 0 , n )$. Put
$U = Z ( \lambda , 0 ) \times Z ( \lambda , 0 )$, then $(x,y) \in
G_{s,0} \cap U$. If $(x',y') \in
G_{s,n} \cap U$ then $x' ( 0 , n ) = \lambda = y' ( 0, n )$ and since
$(x' , y' ) \in G_{s,n}$, we have
$x' ( n ,  \ell ) = y' ( n , \ell )$ for $\ell \ge n$ . Hence
$x' (0, \ell ) = y' ( 0 , \ell )$ for all $\ell \ge 0$ and so
$(x',y') \in G_{s,0} \cap U$, which
shows that $G_{s,0}$ is open in $G_{s,n}$ as required.
% The final statement follows directly from the definitions.
\end{proof}

\begin{rem} \label{opcase}
There is a homeomorphism
$\Lambda^\Delta \to ( \Lambda^{\text{op}} )^\Delta$ given by
$x \mapsto x^{\text{op}}$ where
$$
x^{\text{op}} (m,n) = x ( -n , -m )^{\text{op}} .
$$
Note that for $n \in \Z^k$ and $x \in \Lambda^\Delta$ we have
$
( \sigma^{\text{op}} )^n ( x^{\text{op}} ) = \sigma^{-n} (x)^{\text{op}},
$
where $\sigma^{\text{op}}$ is the
shift action of $\Z^k$ on $( \Lambda^{\text{op}} )^\Delta$.
For every
$x, y \in \Lambda^\Delta$ we have $x \sim_s y$ if and only if
$x^{\text{op}} \sim_u y^{\text{op}}$ and $x \sim_u y$ if and only if
$x^{\text{op}} \sim_s y^{\text{op}}$.
Hence $G_u (\Lambda ) = G_s ( \Lambda^{\text{op}} )$ and
$G_u ( \Lambda^{\text{op}} ) = G_s ( \Lambda )$.
\end{rem}

\begin{rem} \label{asrel}
As in \cite{pt1} (cf.\ \cite{ru2}) we define the asymptotic
relation on $\Lambda^\Delta$ as
follows: For $x,y \in \Lambda$ we put $x \sim_a y$  if $x \sim_s y$
and $x \sim_u y$. Observe that $x \sim_a y$ if and only if there is
$m \in \N^k$ so that for all $n \ge m$ we have
\begin{equation} \label{doubletail}
x ( m , n ) = y ( m , n ) \quad \text{and} \quad
x( - n , -m ) = y ( -n , -m ) .
\end{equation}
Let $G_a$ denote the groupoid derived from this equivalence relation.
We endow it with a topology that makes it an $r$-discrete groupoid.
Given $(x , y ) \in G_a$, there is an $m \in \N^k$ so that
(\ref{doubletail}) holds; set
$\lambda = x ( -m , m )$ and $\nu = y ( -m, m )$. There is a unique map
$\varphi_{\nu , \lambda} : Z ( \lambda , -m ) \to Z ( \nu , -m )$
such that $\varphi_{\nu , \lambda} (x) = y$ and
$$
\varphi_{\nu , \lambda} ( z ) ( m , n )  = z ( m , n )
\quad \text{and} \quad
\varphi_{\nu , \lambda} ( z ) ( - n , -m ) = z ( -n , -m )
$$
for all $z \in Z ( \lambda , - m )$ and $n \ge m$.
Note that the map $\varphi_{\nu , \lambda}$ is a homeomorphism and
$z \sim_a \varphi_{\nu , \lambda} (z)$ for all
$z \in Z ( \lambda , -m )$.
We let $U_{\nu , \lambda} \sub G_a$ be the graph of $\varphi_{\nu , \lambda}$
$$
U_{\nu , \lambda} =
\{ ( \varphi_{\nu , \lambda} (z) , z ) : z \in Z ( \lambda , -m ) \} .
$$
The collection
$\{ U_{\nu , \lambda}\}_{\nu , \lambda}$ forms a
basis for the topology of $G_a$ in which the $U_{\nu , \lambda}$ are
compact open sets. Evidently the restriction of the
range map to each $U_{\nu , \lambda}$ is a homeomorphism onto
$Z ( \nu , -m )$; hence $G_a$ is $r$-discrete.
\end{rem}

\section{Ruelle algebras} \label{ruelle}

As in \cite{pt1} the stable and unstable $C^*$-algebras are given by
$S := C^* ( G_s )$ and $U := C^* ( G_u )$. The Ruelle algebras, $R_s$
and $R_u$, are defined as the crossed products of $S$ and $U$ by the
natural $\Z^k$ actions (see \cite{pts}).  We shall need to
show that $G_s$ and $G_u$ have Haar systems; for this it will be necessary
to invoke a suitable version of the Perron-Frobenius Theorem, for an
irreducible $k$-graph $\Lambda$ with $\Lm^0$ finite (cf.\
\cite{pt2}). As in \cite{pt1} we show that there is a
densely-defined trace on $S$ and $U$ which is scaled by the $\Z^k$
action. Finally, we discuss the corresponding facts in the asymptotic case.

Let $\Lambda$ be a $k$-graph.
For $u, v \in \Lm^0$, $p \in \N^k$ set
$$
%\sideset{_{u}^{}}{_{v}^{n}}{\operatorname*{\Lambda}} =
\Lm^p (u, v) = \{ \lm \in \Lm^p : u = r(\lm) \text{ and }v = s(\lm) \};
$$

then for each $p \in \N^k$ we obtain a nonnegative integer valued
matrix $|\Lm^p|$ indexed by $\Lm^0$ given by
$|\Lm^p|(u, v) = |\Lm^p(u, v)|$ for $u, v \in \Lm^0$.
For $p, q \in \N^k$, we have $|\Lm^{p + q}| = |\Lm^p||\Lm^q|$.
Let $\R_+$ denote the collection of positive real numbers.

\begin{lem}\label{perron}
\textnormal{(cf.\ \cite{pt2})}
Suppose that $\Lm$ is irreducible and $\Lambda^0$ is finite. Then there
exist $t \in \R_+^k$,
$a : \Lm^0 \to \R_+$ and $b : \Lm^0 \to \R_+$ with
$\sum_{v \in \Lm^0} a(v)b(v)=1$ such that
for all $p \in \N^k$ we have
\begin{align}
\sum_{u \in \Lm^0} a(u) |\Lm^p|(u, v) & = t^p a(v) \quad
\text{for all }v \in \Lm^0 \label{pfprop2} \\
\sum_{v \in \Lm^0} |\Lm^p|(u, v)b(v) & = t^pb(u) \quad
\text{for all }u \in \Lm^0 . \label{pfprop1}
\end{align}
\end{lem}
\begin{proof}
Since $\Lambda$ is irreducible, there is an integer matrix $A$ with
all positive entries which may be written as a sum of matrices of the
form $| \Lambda^p |$ for various $p \in \N^k$.
By the Perron-Frobenius theorem (see \cite[Theorem 1.5]{sn} for
example) there are
functions $a,b : \Lm^0 \to \R_+$
satisfying $\sum_{v \in \Lm^0} a(v)b(v)=1$ and a number $T \in \R_+$
such that
\begin{align*}
\sum_{v \in \Lm^0}a(u) A (u, v) &= T a(v)
\text{ for all } v \in \Lm^0 \\
\sum_{u \in \Lm^0} A (u, v)b(v) &= T b(u)
\text{ for all } u \in \Lm^0 .
\end{align*}
For $i=1 , \ldots , k$ let $e_i$
denote the canonical generators of $\N^k$, then since $A$ commutes
with $| \Lambda^{e_i} |$ for each $i$ there exist nonnegative $t_i$
such that the same formulas hold with $A$ replaced by
$| \Lambda^{e_i} |$ and $T$ replaced by $t_i$; formulas
(\ref{pfprop1}) and (\ref{pfprop2}) now follow with
$t = ( t_1 , \ldots , t_k )$.

It remains to show that $t_i >0$ for each $i$. Let $u \in \Lambda^0$,
then by the standing assumption $| \Lambda^e | (u,v) > 0$ for some
$v \in \Lambda^0$; applying  (\ref{pfprop1}) we have
$$
\sum_{v \in \Lm^0} |\Lm^e |(u, v) b(v) = t_1 \cdots t_k b(u) . % \prod_i t_i
$$
Since the left hand side is evidently positive,
$t_1 \cdots t_k >0$; hence $t_i >0$ for all $i$ as required.
\end{proof}

We construct the analog of the Parry measure $\mu$ on $\Lm^\Dl$ as
follows (cf.\ \cite{pt2}).

\begin{pro} \label{mudef}
Suppose that $\Lm$ is irreducible and $\Lambda^0$ is finite. Then
there is a shift invariant probability measure $\mu$ on
$\Lambda^\Delta$ such that
$$
\mu(Z(\lm, n)) = t^{-d(\lm)}a(r(\lm))b(s(\lm)) ,
$$
for all $\lambda \in \Lambda$ and $n \in \Z^k$.
\end{pro}

\begin{proof}
We must show that $\mu$ is well-defined on cylinder sets. Given
$\lambda \in \Lambda$ and $n \in \Z^k$, observe that for $m \ge 0$ we
may write
$Z ( \lambda , n )$ as a disjoint union by expanding on the right:
$$
Z ( \lambda , n ) =
\coprod_{\substack{ \nu \in \Lambda^{m} \\ r ( \nu ) = s ( \lambda )}}
Z ( \lambda \nu , n ) .
$$
Then we compute $\mu$ of the right-hand side (using (\ref{pfprop1}))
\begin{align*}
\sum_{\substack{ \nu \in \Lambda^{m} \\ r ( \nu ) = s ( \lambda )}}
\mu ( Z ( \lambda \nu , n ) ) &=
\sum_{\substack{ \nu \in \Lambda^{m} \\ r ( \nu ) = s ( \lambda )}}
t^{-d ( \lambda \nu )}a ( r ( \lambda ) ) b ( s ( \nu ) ) \\
 &= t^{-d ( \lambda)} a ( r ( \lambda ) ) \sum_{v \in \Lm^0}
t^{-m} \sum_{\nu \in \Lambda^{m}(s(\lm), v )} b ( s ( \nu ) )\\
 &= t^{-d ( \lambda)} a ( r ( \lambda ) ) \sum_{v \in \Lm^0}
 t^{-m}|\Lambda^{m}|(s(\lm), v ) b ( v )\\
 &=  t^{-d ( \lambda)} a ( r ( \lambda ) )  b ( s ( \lm) )
 = \mu(Z(\lm, n)).
\end{align*}
If we write $Z ( \lambda , n )$ as a disjoint union by expanding on
the left:
$$
Z ( \lambda , n ) =
\coprod_{\substack{ \nu \in \Lambda^{m} \\ s ( \nu ) = r ( \lambda )}}
Z ( \nu \lambda , n-m ) ,
$$
then a similar calculation (using
(\ref{pfprop2})) shows that $\mu$ of each side is the same and
completes the demonstration that $\mu$ is well-defined.
Thus $\mu$ extends to a probability measure which is invariant
under the action of $\Z^k$. % by construction.
\end{proof}

Our next task is to decompose $\mu$ locally into measures which give
rise to Haar systems for $G_s$ and $G_u$ .

Fix $x \in \Lambda^\Delta$; then for %$m \ge 0$  and
$\lambda \in \Lambda$ with $s(\lm) = x(0)$
we define
$$
Z^- (\lambda , x ) = E_x \cap Z (\lambda, -d(\lambda) ) .
$$

\noindent
Likewise for % $n \ge 0$ and $\lambda \in \Lambda^n$
$\lambda \in \Lambda$ with $r(\lm) = x(0)$
we define
$$
Z^+ ( \lambda , x ) = F_x \cap Z (\lambda, 0 ) .
$$

\begin{pro}\label{mus}
\textnormal{(cf.\ \cite{pt2})}
Suppose that $\Lm$ is irreducible and $\Lambda^0$ is finite. Then for
each $x \in \Lambda^\Delta$ there is a measure $\mu_s^x$ on $E_x$
and a measure $\mu_u^x$ on $F_x$ such that
$$
\mu_s^x ( Z^- ( \lambda , x ) ) = t^{-d ( \lambda )} a (r(\lambda))
\quad \text{and} \quad \mu_u^x ( Z^+ ( \lambda , x ) ) = t^{-d (
\lambda )} b ( s ( \lambda ) )
$$
for all $\lambda \in \Lambda$. The restriction of $\mu$ to $Z ( x(0), 0 )$ is
$\mu_s^x \times \mu_u^x$ after identifying
$E_x \times F_x$ with $Z ( x(0), 0 )$ as in Proposition \ref{bracket}.
Moreover for $p \in \N^k$ we have
\begin{equation}\label{trans}
\mu_s^x = t^p\mu_s^{\sigma^p x} \circ \sigma^p
\quad \mbox{and} \quad
\mu_u^x = t^p \mu_u^{\sigma^{-p} x}\circ \sigma^{-p}
\end{equation}
on $E_x$ and $F_x$ respectively.
\end{pro}

\begin{proof}
The first assertion is clear.
% Every cylinder set contained in $Z (x(0),0)$ is
For the next part it suffices to consider cylinder sets
of the form $Z ( \lambda , n )$ where
$\lambda = \lambda^- \lambda^+ \in \Lambda$ with
$r(\lambda^+) = x(0)$ and $n = - d ( \lambda^- )$.  After identifying
$Z^+ ( \lambda^+ , x ) \times Z^- ( \lambda^- , x )$ with
$Z ( \lambda , - d ( \lambda^- ))$ (as in Proposition \ref{bracket}) we have
\begin{align*}
\mu_s^x \times \mu_u^x ( Z ( \lambda , - d ( \lambda^- )) ) &=
\mu_s^x( Z^- ( \lambda^- , x ) )\mu_u^x ( Z^+ ( \lambda^+ , x ) )\\
&= t^{-d (\lambda^-)} a (r(\lambda^-))t^{-d (\lambda^+)} b (s (\lambda^+ ) )\\
&= t^{ -d (\lambda )} a ( r ( \lambda ) ) b ( s ( \lambda ) )  \\
&= \mu ( Z (\lambda, -d(\lambda^-)) );
\end{align*}

\noindent
hence the restriction of $\mu$ to $Z (x(0), 0)$ is
$\mu_s^x \times \mu_u^x$ as required.  From the definitions it is
straightforward to verify that for $p \in \N^k$
$$
\sigma^p Z^- ( \lambda , x ) = Z^- ( \lambda x (0,p) , \sigma^p x )
\quad \text{and} \quad
\sigma^{-p} Z^+ ( \lambda , x ) = Z^+ ( x (-p,0) \lambda , \sigma^{-p}
x ) .
$$
Hence
\begin{align*}
\mu_s^{\sigma^p x} ( \sigma^p Z^- ( \lambda , x ) ) &= t^{-p} \mu_s^x ( Z^- ( \lambda , x ) )  \\
\mu_u^{\sigma^{-p} x} ( \sigma^{-p} Z^+ ( \lambda , x ) ) &= t^{-p}
\mu_u^x ( Z^+ ( \lambda , x ) ) ;
\end{align*}
equations (\ref{trans}) then follow on $E_x$ and $F_x$ respectively.
\end{proof}

Note that for $x \in \Lambda^\Delta$ we have that
$E_x = \{ y : \pi (y) = \pi (x) \}$; evidently if $\pi (y) = \pi (x)$ then
$E_y = E_x$ and $\mu_s^y = \mu_s^x$. For $z = \pi (x)$ let
$\mu_{s,0}^z$ denote the extension of $\mu_s^x$ to $\Lambda^\Delta$
($\mu_{u,0}^z$ is defined similarly). Observe that
$\mu_{s,0} = \{ \mu_{s,0}^z : z \in \Lambda^\Omega \}$ is a full
$\pi$-system: The
continuity of the system follows from the fact that
$$
\mu_{s,0} \left( \chi_{Z ( \lambda , n )} \right)(z) = t^{-d ( \lambda )}
a ( r ( \lambda ) )
$$
is locally constant.
By Propositions \ref{pisys} and \ref{gsprops} we obtain a Haar system
$$
\tilde{\mu}_{s,0} = \{ \tilde{\mu}_{s,0}^x : x \in \Lambda^\Delta \}
$$
for $G_{s,0}$. Recall from Proposition \ref{gsprops} that for
$p \ge 0$ we have
$G_{s,p} = \Lambda^\Delta \star_{\pi \circ \sigma^p} \Lambda^\Delta$.

\begin{pro} \label{haarhaar}
Let $\Lambda$ be an irreducible $k$-graph with $\Lambda^0$ finite.
For $p \ge 0$ and $x, y \in \Lambda^\Delta$ if
$\pi \circ \sigma^p x = \pi \circ \sigma^p y$ then
\begin{equation} \label{pcase}
\sigma^{-p} E_{\sigma^p x} = \sigma^{-p} E_{\sigma^p y} \quad
\text{and} \quad
\mu_s^{\sigma^p x} \circ \sigma^p =
\mu_s^{\sigma^p y} \circ \sigma^p .
\end{equation}
For
$z = \pi\circ \sigma^p (x)$ let $\mu_{s,p}^z$ denote the
extension of the measure $t^p \mu_s^{\sigma^p x} \circ \sigma^p$ from
$\sigma^{-p} E_{\sigma^p x}$ to $\Lambda^\Delta$.  Then
$\mu_{s,p} = \{ \mu_{s,p}^z : z \in \Lambda^\Omega \}$ is a full
$\pi \circ \sigma^p$-system and we obtain a Haar system
$\tilde{\mu}_{s,p} = \{\tilde{\mu}_{s,p}^x : x \in \Lambda^\Delta \}$
for $G_{s,p}$.
If $0 \le p \le q$, the restriction of $\tilde{\mu}_{s,q}^x$ to
$G_{s,p}$ agrees with $\tilde{\mu}_{s,p}^x$ for all
$x \in  \Lambda^\Delta$.  We obtain thereby a Haar system for $G_s$,
$\tilde{\mu}_{s} = \{\tilde{\mu}_{s}^x : x \in \Lambda^\Delta \}$.
A Haar system for $G_u$ is obtained in a similar way.
\end{pro}

\begin{proof}
Note that for $p \ge 0$ and $x, y \in \Lambda^\Delta$,
$y \in \sigma^{-p} E_{\sigma^p x}$ if and only if
$\pi \circ \sigma^p x = \pi \circ \sigma^p y$;
hence, formulas (\ref{pcase}) hold.
By arguments similar to those given immediately prior to the statement
of this proposition, $\mu_{s,p}$ is a full $\pi \circ \sigma^p$-system
and $\tilde{\mu}_{s,p}$ is a Haar system for $G_{s,p}$. The
compatibility of the Haar systems now follows from
a short calculation involving (\ref{trans}). The Haar systems
$\mu_{s,p}$ for $p \ge 0$ may therefore be
patched together to give a Haar system for $G_s$:
$\tilde{\mu}_{s} = \{\tilde{\mu}_{s}^x : x \in \Lambda^\Delta \}$.
\end{proof}

% We now extend these measures to the equivalence classes of $x$
% in the stable and unstable equivalence relations. One checks that that
% this definition is independent of the choice of $x$ in the
% class. Whence we shall obtain Haar systems for the groupoids $G_s$,
% $G_u$.
% Observe that
% $$
% [ x ]_s = \bigcup_{p \ge 0} \sigma^{-p} E_{\sigma^p x} \quad
% \mbox{and} \quad
% [ x ]_u = \bigcup_{p \ge 0} \sigma^p F_{\sigma^{-p} x} ;
% $$

% for $p \ge 0$ by (\ref{trans}) we extend $\mu_s^x$ to
% $\sigma^{-p} E_{\sigma^p x}$
% and extend $\mu_u^x$ to $\sigma^{-p} F_{\sigma^{-p} x}$.
% It also follows that these extensions are compatible.
% If $x \sim_s y$ then there is $N \ge 0$ so that $\sigma^p (x) =
% \sigma^p (y)$ for all $p \ge N$ and then by the above we have
% $\mu_s^x = \mu_s^y$. A similar argument applies to the unstable case.
% We have now all but proved the following:% Proposition:

Now we may define the stable and unstable algebras associated to an
irreducible $k$-graph $\Lambda$ with finitely many vertices
$$
S = S ( \Lambda ) = C^* ( G_s ( \Lambda ) ) \quad \mbox{and} \quad
U = U ( \Lambda ) = C^* ( G_u ( \Lambda ) ) .
$$
For $n \in \Z^k$ the map $\sigma^n \times \sigma^n$ yields an
automorphism of the
stable and unstable equivalence relations but it rescales the Haar
systems by (\ref{trans}); indeed
\begin{equation} \label{actionscale}
\tilde{\mu}_s^x \circ ( \sigma^{-n} \times \sigma^{-n} ) =
t^n  \tilde{\mu}_s^{\sigma^n x} \quad \text{and}  \quad
\tilde{\mu}_u^x \circ ( \sigma^{-n} \times \sigma^{-n} ) =
t^{-n} \tilde{\mu}_u^{\sigma^n x} .
\end{equation}
This induces actions $\beta_s , \beta_u$
of $\Z^k$ on both $S$ and $U$ given for $n \in \Z^k$ by
\begin{align*}
\beta^n_s ( f ) (x,y) &= t^{n} f ( \sigma^{-n} (x) , \sigma^{-n}
(y)) \mbox{ where }  f \in C_c ( G_s )  , \\
\beta^n_u ( f ) (x,y) &= t^{-n} f ( \sigma^{-n} (x) , \sigma^{-n}
(y)) \mbox{ where } f \in C_c ( G_u ) ,
\end{align*}
and extending by continuity to the completions.

The measure $\mu$ on $\Lambda^\Delta$ gives rise to a densely-defined
trace on $S$ as follows:

\begin{pro} \label{mutrace}
Let $\Lambda$ be an irreducible $k$-graph with $\Lambda^0$ finite,
$\mu$ the Parry measure, and
$G_s$ the stable groupoid. For $f \in C_c ( G_s )$ set
$$
\tau_s (f) = \int_{\Lambda^\Delta} f(x,x) d \mu (x) .
$$
Then $\tau_s$ is a densely-defined trace on $C^* ( G_s )$.
A densely-defined trace $\tau_u$ on $C^* ( G_u )$ is defined
similarly. Moreover,  for $n \in \Z^k$ we have
\begin{equation} \label{trscale}
\tau_s \circ \beta_s^n = t^n \tau_s \quad \text{and} \quad \tau_s \circ
\beta_u^n = t^{-n} \tau_u .
\end{equation}
\end{pro}

\begin{proof}
The formulas (\ref{trscale}) follow from (\ref{actionscale}).
We show that $\tau_s$ is a densely-defined trace; the case of $\tau_u$
is similar. It suffices
to show that $\tau_s$ satisfies the trace property ($\tau_s$ is
clearly densely-defined, linear and positive).
For $f, g \in C_c ( G_s )$
there is $p \ge 0$ such that the supports of $f$ and
$g$ are contained in $G_{s,p}$ (by Proposition \ref{gsprops},
$\{ G_{s,p} : p \ge 0 \}$ forms an open cover of $G_s$). Since
$\mu$ decomposes locally as a
product measure as in Proposition \ref{mus} with
$t^p \mu_s^{\sigma^p x} \circ \sigma^p$ in place of $\mu_s^x$,
there is a measure
$\nu_{s,p}$ on $\Lambda^\Omega$ such that
$$
\int_{\Lambda^\Delta} h (x) d \mu (x) =
% \int_{\Lambda^\Omega} \int h ( y ) d \mu_{s,p}^z (y) d \nu_{s,p} (z) .
\int_{\Lambda^\Omega} \mu_{s,p} (h) d \nu_{s,p}
$$
for all $h \in C ( \Lambda^\Delta )$. It follows that
\begin{align*}
\tau_s (fg) &= \int_{\Lambda^\Delta}( fg) (x,x) d \mu (x)
= \int_{\Lambda^\Omega} \int (fg) ( y,y ) d \mu_{s,p}^z (y)
d\nu_{s,p}(z) \\
&= \int_{\Lambda^\Omega} \tau_{\mu_{s,p}} (fg) d \nu_{s,p}
= \int_{\Lambda^\Omega} \tau_{\mu_{s,p}} (gf) d \nu_{s,p}
= \tau_s (gf) ,
\end{align*}
where the penultimate equality follows from Proposition \ref{pisys}.
% Let $f, g \in C_c ( G_s )$,
% then since the $G_{s,m}$ form an
% increasing family of open sets whose union is $G_s$ (see Proposition
% \ref{gsprops}) we have that there
% is some $m$ such that $f , g \in C_c ( G_{s,m} )$.
% We may, without loss of generality assume that $m=0$, for
% $G_{s,m} \cong G_{s,0}$ (by means of an automorphism of $G_s$).
\end{proof}

Let $\Lambda$ be an irreducible $k$-graph with finitely many vertices.
The Ruelle algebras associated to $\Lambda$ are defined
to be the corresponding crossed products (cf.\ \cite{pts, pt2})
$$
R_s = R_s ( \Lambda ) = S ( \Lambda )\times_{\beta_s} \Z^k \quad
\mbox{ and } \quad
R_u = R_u ( \Lambda ) = U ( \Lambda ) \times_{\beta_u} \Z^k .
$$

We express the Ruelle algebras as $C^*$-algebras of the semidirect
product groupoids $G_s \times \Z^k$ and $G_u \times \Z^k$. The unit
space is identified with $\Lambda^\Delta$ via the map
$x \mapsto ((x,x),0)$. The  structure maps are given by
% \begin{align*}
% r ( (x,y) , n ) = ( x , 0 ) ,\quad &
% s ( (x,y) , n ) = ( \sigma^{n} y , 0 ) ,&
% ((x,y),n)^{-1} = ((\sigma^n y , \sigma^n x),-n),  \\
% ( (x,y) , n ) ( ( \sigma^n y , \sigma^n z ) , m ) & = ((x,z),n+m).
% \end{align*}
$$
r ( (x,y) , n ) = x  ,\quad
s ( (x,y) , n ) = \sigma^{n} y , \quad
( (x,y) , n ) ( ( \sigma^n y , \sigma^n z ) , m )  = ((x,z),n+m).
$$
The structure maps of $G_u \times \Z^k$ are defined similarly.
% and $((x,y),n)^{-1} = ((\sigma^n y , \sigma^n x),-n)$.
% We now examine the structure of $G_s \times_{\beta_s} \Z^k$ more
% closely:

\begin{lem} \label{haarhaarsemi}
If $\Lambda$ is an irreducible $k$-graph with $\Lambda^0$ finite, then
both $G_s \times \Z^k$ and $G_u \times \Z^k$
have Haar systems. Moreover, we have
$R_s \cong C^* ( G_s \times \Z^k )$ and
$R_u \cong C^* ( G_u \times \Z^k )$.
\end{lem}

\begin{proof}
Let $\vartheta$ be the measure on $\Z^k$ given by
$\vartheta ( \{ n \} ) = t^{-n}$; then a direct computation using
(\ref{actionscale}) shows that
$\{ \tilde{\mu}_s^x \times \vartheta : x \in \Lambda^\Delta \}$ is a
Haar system for $G_s \times \Z^k$.
\end{proof}

\begin{rem} \label{op}
Let $\Lambda$ be an irreducible $k$-graph with $\Lambda^0$ finite,
then by Remark \ref{opcase} we have
$G_u ( \Lambda ) = G_s ( \Lambda^{\text{op}} )$ and
$G_s ( \Lambda ) = G_u ( \Lambda^{\text{op}} )$.
Note that $\Lambda^{\text{op}}$ is also irreducible and
$( \Lambda^{\text{op}} )^0 = \Lambda^0$ is finite.
We have
$U ( \Lambda ) = S ( \Lambda^{\text{op}} )$ and
$S ( \Lambda ) = U ( \Lambda^{\text{op}} )$, similarly
$R_u ( \Lambda ) = R_s ( \Lambda^{\text{op}} )$ and
$R_s ( \Lambda ) = R_u ( \Lambda^{\text{op}} )$.
Henceforth we shall focus our attention on the stable case.
\end{rem}

\begin{rem} \label{asruelle}
The asymptotic $C^*$-algebra may also be defined:
$$
A = A ( \Lambda )=  C^* ( G_a ( \Lambda ) ) .
$$
Note that since $G_a$ is $r$-discrete, it has a Haar system consisting
of counting measures. The asymptotic Ruelle algebra is defined
as the crossed product
$$
R_a = R_a ( \Lambda ) = A ( \Lambda ) \times_{\beta_a} \Z^k
= C^* ( G_a ( \Lambda ) \times \Z^k ) .
$$
Suppose that $\Lambda$ is irreducible and $\Lambda^0$ is
finite. With notation as in Remark \ref{asrel},
$$
\mu ( Z ( \lambda , -m ) ) =
t^{-2m} a ( r ( \lambda )) b ( s ( \lambda ) ) =
t^{-2m} a ( r ( \nu )) b ( s ( \nu ) ) =
\mu ( Z ( \nu , -m ) ) ;
$$
hence $\mu$ is invariant under $G_a$. Thus we may
define a unital trace on $A$ by
$$
\tau_a (f) = \int f (x,x) d \mu (x)
$$
for $f \in C_c ( G_a )$. The $\Z^k$ action $\sigma$ on
$\Lambda^\Delta$ induces an action $\sigma \times \sigma$ on $G_a$ and
hence we get an action  $\beta_a : \Z^k \to \Aut A$.
Since $\tau_a$ is invariant under $\beta_a$, we also obtain a trace on
$R_a$.
\end{rem}

\section{Morita Equivalence}\label{morita}

In this section we prove our main results. If $\Lambda$ is irreducible
and $\Lambda^0$ is finite
we show that the stable algebra $S$ is strongly
Morita equivalent to $C^* ( \Lambda )^\alpha$ and that $R_s$ is
strongly Morita equivalent to $C^* ( \Lambda )$.  Hence, if $\Lambda$
satisfies the aperiodicity condition, then $R_s$ is a stable Kirchberg
algebra.
% This will allow us to deduce various properties of
% the Ruelle algebras from the corresponding properties of
% $C^* ( \Lambda )$.
% We then conclude by discussing the asymptotic case; we
% show (cf.\ \cite{pt1, pt2}) that the asymptotic algebra $A$ is strongly
% Morita equivalent to $S \otimes U$.

We begin by stating a groupoid equivalence result that will be useful in both
cases. This result is no doubt well known to the experts but, as we are
unable to find an explicit reference, we provide a proof.

Let $\Gamma$ be a locally compact Hausdorff groupoid. Given a right
principal $\Gamma$-space $Y$, one may construct the imprimitivity
groupoid $Y \star_\Gamma Y^{\text{op}}$ where $Y^{\text{op}}$ is the
corresponding left principal $\Gamma$-space. By
\cite[Theorem 3.5]{mw2} (see also \cite[Theorem 5.31]{m}) $Y$
implements an equivalence between the imprimitivity groupoid and
$\Gamma$ in the sense of \cite[Definition 2.1]{mrw}.

Let $X$ be a locally compact Hausdorff space and let
$\psi : X \rightarrow \Gamma^0$ be a
continuous open surjection. Set
\begin{equation} \label{z}
Z = X \star\Gamma = \{ ( x , \gamma ) : x \in X , \gamma \in \Gamma ,
\psi (x) = r ( \gamma ) \} ;
\end{equation}
we define a right action of $\Gamma$ on $Z$ as follows:
$s : Z \to \Gamma^0$ is given by $s ( x , \gamma ) = s ( \gamma )$ and
the map $Z \star \Gamma \to Z$ by
$( (x ,\gamma_1 ) , \gamma_2 ) = ( x , \gamma_1 \gamma_2 )$. There is
a corresponding left action of
$\Gamma$ on $Z^{\text{op}} = \Gamma \star X$.

\begin{lem} \label{imp}
With the above structure maps $Z$ is a right principal $\Gamma$-space.
Moreover the imprimitivity groupoid $Z \star_\Gamma Z^{\text{op}}$ is
isomorphic to
$$
\Gamma^\psi := \{ (x,\gamma,y) : x,y \in X , \gamma \in \Gamma ,
\psi (x) = r ( \gamma ) , \psi (y) = s ( \gamma ) \} ,
$$
equipped with the relative topology.
Therefore $Z$ implements an equivalence between the groupoids $\Gamma$
and $\Gamma^\psi$.
\end{lem}

\begin{proof}
To show that $Z$ is a right principal $\Gamma$-space we must show that
the action is free and proper. The action is clearly free (because the
action of a groupoid on itself is free). It suffices to show that the
map  $Z \star \Gamma \rightarrow Z \times Z$ given by
$$ ( (  x , \gamma_1 ) , \gamma_2 )
\mapsto ( ( x , \gamma_1 \gamma_2 ) , ( x, \gamma_1 ) )$$
is a
homeomorphism onto a closed set (see \cite[Lemma 2.2]{mw2}).
This follows from a
similar fact for the right action of a groupoid on itself.
By \cite[Theorem 3.5]{mw2} $Z$ is a
groupoid equivalence between the imprimitivity groupoid
$Z \star_\Gamma Z^{\text{op}}$ and $\Gamma$.
The isomorphism from  $Z \star_\Gamma Z^{\text{op}}$ to $\Gamma^\psi$
is given by the map
$$
( (x , \gamma_1 ) , ( \gamma_2 , y ) ) \mapsto
( x , \gamma_1 \gamma_2 , y ).
$$
The result now follows from this identification.
\end{proof}

The construction of $\Gamma^\psi$ above appears in
\cite[Proposition 5.7]{ku1} (though in a more specialized setting).
Recall the the restriction map $\pi : \Lambda^\Delta \to
\Lambda^\Omega$ is a continuous open surjection.

\begin{lem} \label{gsisgood}
Let $\Lambda$ be a $k$-graph and $G_s$ be the stable groupoid associated
to $\Lambda$. Then the map
$(x,y) \mapsto ( x , ( \pi (x) , 0 , \pi (y) ) , y )$ gives an isomorphism
$G_s \cong \left( \Gamma_\Lambda \right)^\pi$.
\end{lem}

\begin{proof}
Recall that $(x,y) \in G_{s,m}$ precisely when
$\pi ( \sigma^m x)= \pi ( \sigma^m y)$. For $m \ge 0$ the given map
restricts to a
homeomorphism from $G_{s,m}$ to $( \Gamma_{\Lambda,m} )^\pi$ where
$\Gamma_{\Lambda,m}$ is the subgroupoid of $\Gamma_\Lambda$ formed by
those $(x,0,y)$ where $\sigma^m (x) = \sigma^m (y)$. Since the
topology on $\Gamma_\Lambda$ is equivalent to the inductive limit
topology from these subgroupoids, the given map is a homeomorphism.
It is routine to check that the map is a groupoid morphism.
\end{proof}

The groupoid equivalence (in the sense of \cite{mrw}) between $G_s$
and $\Gamma_\Lambda$ now follows.

% $$
% Z = \Lambda^\Delta \star \Gamma_\Lambda = \{ ( x , (y,z)) : x \in
% \Lambda^\Delta , (y,z) \in \Gamma_\Lambda , \pi (x) = y  \}
% $$
% and give $Z$ the relative topology.

\begin{thm} \label{zisequiv}
The space $Z = \Lambda^\Delta \star \Gamma_\Lambda$  is a $( G_s ,
\Gamma_\Lambda )$-equivalence.
In particular $G_s$ is amenable in the sense of \cite{ar}.
Moreover if $\Lambda$ is irreducible and $\Lambda^0$ is finite then
the stable algebra $S$ is strongly Morita equivalent to
$C^* ( \Lambda )^\alpha$, and therefore is an AF-algebra.
\end{thm}

\begin{proof}
The first part follows from Lemmas \ref{imp},  \ref{gsisgood}. That
$G_s$ is amenable now follows from \cite[Theorem 2.2.17]{ar}.
If $\Lambda$ is irreducible and $\Lambda^0$ is finite then by
Proposition \ref{haarhaar}, $G_s$ has a Haar system ($\Gamma_\Lambda$
has a Haar system consisting of counting measures) so that by
\cite[Theorem 2.8]{mrw} $S$ is strongly Morita equivalent
to $C^* ( \Gamma_\Lambda ) = C^* ( \Lambda )^\alpha$. The final
assertion then follows from \cite[Lemma 3.2]{kp}.
\end{proof}

We could have deduced that $S = C^* ( G_s )$ is AF more directly: It
follows from the fact that
$$
C^* ( G_s ) = \lim_m C^* ( G_{s,m} )
$$
and that $C^* ( G_{s,m} )$ is strongly Morita equivalent to the abelian AF
algebra $C ( \Lambda^\Omega )$ for each $m$ (see Propositions \ref{pisys}
and \ref{gsprops}).
If $\Lambda$ is primitive then we can say more.

\begin{cor} \label{prim}
Let $\Lambda$ be a primitive $k$-graph, then
$C^* ( \Lambda )^\alpha$ is a simple AF algebra and hence so is $S$.
\end{cor}

\begin{proof}
Suppose that $\Lambda$ is primitive; then by \ref{star},
$\Lambda^0$ is finite. Moreover, there is an $n >0$ (i.e.\ all
coordinates are positive) such that for every
$u, v \in \Lambda^0$ there is $\lambda \in \Lambda^n$ with $u =
s(\lambda )$ and $v = r ( \lambda )$. It follows that all the entries
of the matrix $| \Lambda^n |$ are positive. Since the sequence
$\{ j n : j \in \N \}$ is cofinal in $\N^k$, we have that
$C^* ( \Lambda )^\alpha = \lim_{j}
\mathcal{F}_{jn}$.
The multiplicity
matrix of the inclusions may be identified with $| \Lambda^n |$ and
the result now follows from \cite[Corollary 3.5]{b}.
\end{proof}
Analogous assertions hold for $U$ when $\Lambda$ is replaced
by $\Lambda^{\text{op}}$ (see Remark \ref{op}).
%The case of the stable Ruelle algebras is similar:

\begin{lem} \label{gsxzkisgood}
Let $\Lambda$ be a $k$-graph, $G_s$ be the associated stable groupoid
and $G_s \times \Z^k$ be the semidirect product groupoid (see Lemma
\ref{haarhaarsemi}). Then the map
$$
\varphi : ((x,y),n) \mapsto ( x , ( \pi (x) , n , \pi ( \sigma^{n} y) ) ,
\sigma^{n} y )
$$
gives an isomorphism
$G_s \times \Z^k \cong \left( \mathcal{G}_\Lambda \right)^\pi$.
\end{lem}

\begin{proof}
That $( \pi (x) , n , \pi ( \sigma^n y ) ) \in \mathcal{G}_\Lambda$
follows from (\ref{stationary}); so $\varphi$ is well-defined and
evidently injective. Given
$(x,(\pi (x) , n , \pi (z) ) , z ) \in ( \mathcal{G}_\Lambda )^\pi$
we have
$$
(x,(\pi (x) , n , \pi (z) ) , z ) = \varphi ((x,\sigma^{-n} z),n) ;
$$
hence $\varphi$ is surjective.
Recall that $G_s \times \Z^k$ is given the product topology
and observe that the restriction of $\varphi$ to $G_s \times \{ 0 \}$
agrees with the homeomorphism defined in Lemma
\ref{gsisgood}. Similarly the restriction to $G_s \times \{ n \}$ is a
homeomorphism onto the set
$$
\{ (x,(\pi(x),n,\pi(y)),y) : x ,y \in \Lambda^\Delta,\sigma^\ell \pi (x)
= \sigma^m \pi (y) , n = \ell - m \}.
$$
The reader is invited to check that the map is a groupoid morphism.
\end{proof}

% For $v \in \Lambda^0$ fix $x_v \in \Lambda^\Delta$ with $x(0)=v$. Set
% $$
Recall that a $k$-graph is said to satisfy the {\em aperiodicity condition}
if for every vertex $v$ there is $x \in \Lambda^\Omega$ with $x(0)=v$
which is not eventually periodic (see \cite[Definition 4.1]{kp}).
Let $\mathcal{N}$ denote the bootstrap class of $C^*$-algebras
to which the UCT applies (see \cite{rsc}).

\begin{thm} \label{zequiv2}
The space $Z = \Lambda^\Delta \star \mathcal{G}_\Lambda$  is a
$( G_s \times \Z^k , \mathcal{G}_\Lambda )$-equivalence.
Therefore if $\Lambda$ is irreducible and $\Lambda^0$ is finite, then
the stable Ruelle algebra $R_s$ is strongly Morita equivalent to
$C^* ( \Lambda )$ and hence is nuclear and lies in the bootstrap class
$\mathcal{N}$. If in addition $\Lambda$ satisfies the
aperiodicity condition, then $R_s$ is simple, stable and purely
infinite.
Hence the isomorphism class of $R_s$ is
completely determined by $K_* ( R_s ) = K_* ( C^* ( \Lambda ) )$.
\end{thm}
\begin{proof}
The first part follows from Lemmas \ref{imp},  \ref{gsxzkisgood}.
If $\Lambda$ is irreducible and $\Lambda^0$ is finite, then by Lemma
\ref{haarhaarsemi}, $G_s \times \Z^k$ has a Haar system
($\mathcal{G}_\Lambda$ has a Haar system consisting of counting measures)
so that by \cite[Theorem 2.8]{mrw} $R_s$ is strongly Morita equivalent
to $C^* ( \mathcal{G}_\Lambda ) = C^* ( \Lambda )$.
By \cite[Theorem 5.5]{kp} $R_s$ is nuclear and lies in the bootstrap
class $\mathcal{N}$ (since strong Morita equivalence preserves these
properties). If $\Lambda$
is irreducible then it is clearly cofinal and if in addition $\Lambda$
satisfies the aperiodicity condition, it follows from
\cite[Proposition 4.8]{kp} that $C^* ( \Lambda )$
is simple. For every vertex $v \in \Lambda$ there is a morphism
$\lambda \in \Lambda$ with $d ( \lambda ) \ne 0$ such that
$r ( \lambda ) = s ( \lambda ) = v$ and so $C^* ( \Lambda )$ is purely
infinite by
\cite[Proposition 4.9]{kp}. By Zhang's dichotomy a simple purely
infinite $C^*$-algebra is either unital or stable (see
\cite[Theorem 1.2]{z}); since $R_s$ is not unital it must be stable.
The Kirchberg-Phillips Theorem applies and the isomorphism
class of $R_s$ is completely determined by $K_* ( R_s )$
(see \cite[Theorem C]{ki}, \cite[Corollary 4.2.2]{ph}).
\end{proof}
An analogous assertion holds for $R_u$ when $\Lambda$ is replaced by
$\Lambda^{\text{op}}$.
The aperiodicity condition is necessary in the statement of the above theorem.
There is an irreducible $2$-graph $\Lambda$ with one vertex  which is
not aperiodic --- every path has period $(1,-1)$. Furthermore,
$C^* ( \Lm  ) \cong \mathcal{O}_2 \otimes C ( \T )$
is neither simple nor purely infinite (see
\cite[Example 6.1]{kp}).

The restriction of Theorem \ref{zequiv2} to the case $k=1$ is
certainly well known, but we have been unable to find a reference.

\begin{cor}
Let $A \in M_n ( \N )$ be irreducible and $R_s$ be the
stable Ruelle algebra of the associated Markov shift. Then
$R_s$ is strongly Morita equivalent to $\mathcal{O}_A$.
\end{cor}

\begin{rem}
Suppose that $\Lambda$ is an irreducible $k$-graph with $\Lambda^0$
finite. Then
the $2k$-graph $\Lambda \times \Lambda^{\text{op}}$ is irreducible and
$( \Lambda \times \Lambda^{\text{op}} )^0 =
\Lambda^0 \times ( \Lambda^{\text{op}} )^0$ is finite. We have
$( \Lambda \times \Lambda^{\text{op}} )^\Delta =
\Lambda^\Delta \times ( \Lambda^{\text{op}} )^\Delta$ and
$$
G_s ( \Lambda \times \Lambda^{\text{op}} ) =
G_s ( \Lambda ) \times G_s ( \Lambda^{\text{op}} ) \cong
G_s ( \Lambda ) \times G_u ( \Lambda ) ;
$$
hence
$$
S ( \Lambda \times \Lambda^{\text{op}} )
\cong S ( \Lambda ) \otimes U ( \Lambda ) .
$$
Moreover $A ( \Lambda )$ is strongly Morita equivalent to
$S ( \Lambda ) \otimes U ( \Lambda )$ as in \cite[Theorem 3.1]{pt1}.
The ``same'' argument applies: Define a map
$\phi : \Lambda^\Delta \to \Lambda^\Delta \times
( \Lambda^{\text{op}} )^\Delta$ by
$x \mapsto ( x , x^{\text{op}} )$, then $N = \phi ( \Lambda^\Delta )$
is an abstract transversal of the groupoid
$G_s ( \Lambda \times \Lambda^{\text{op}} )$ in the sense of
\cite[Example 2.7]{mrw}. Furthermore, $G_a ( \Lambda )$ is isomorphic
to the reduction $G_s ( \Lambda \times \Lambda^{\text{op}} )^N_N$ (for
$x \sim_a y$ if and only if $x \sim_s y$ and
$x^{\text{op}} \sim_s y^{\text{op}}$, that is,
$(x,x^{\text{op}} ) \sim_s ( y , y^{\text{op}} )$). It follows that
$A ( \Lambda )$ is an AF algebra and if $\Lambda$ is primitive then
$A ( \Lambda )$ is simple. However, $R_a ( \Lambda )$ is not purely
infinite since it has a trace (see Remark \ref{asruelle}).
\end{rem}


\begin{thebibliography}{KMRW}\label{biblio}
\bibitem[AR]{ar}
C. Anantharaman-Delaroche and J. Renault, {\em Amenable groupoids.}
%With a foreword by Georges Skandalis and Appendix B by E. Germain.
Monographies de L'Enseignement Math\'{e}matique {\bf 36},
L'Enseignement Math\'{e}matique, Geneva, 2000.
% Amenability is invariant under equivalence - Theorem 2.2.17
% \bibitem[BPRS]{bprs}
% Teresa Bates, David Pask, Iain Raeburn, and Wojciech Szyma\'{n}ski
% {\em The $C^*$-Algebras of Row-Finite Graphs,}
% New York Journal of Mathematics,  {\bf 6} (2000) 307-324.
\bibitem[B]{b}
O. Bratteli, {\em Inductive limits of finite dimensional
$C^*$-algebras}, Trans. Amer. Math. Soc. {\bf 171} (1972) 195-234.
\bibitem[CK]{ck}
J. Cuntz and W. Krieger, {\em A class of $C^*$-algebras and
topological Markov  chains,} Invent. Math. {\bf 56} (1980) 251-268.
\bibitem[KPS]{kps}
J. Kaminker, I. Putnam and J. Spielberg, {\em Operator algebras and
hyperbolic dynamics,} in Operator Algebras and Quantum Field Theory
(Accademia Nazionale dei Lincei, Roma, %  July 1-6
1996) S. Doplicher et al., %J. Roberts, L. Zsido,
Eds., 525-532, International Press 1997.
\bibitem[D]{d}
V. Deaconu,
{\em $C^*$-algebras associated to higher dimensional Markov shifts},
in preparation.
% \bibitem[D1]{d1}
% V. Deaconu, {\em Groupoids associated with endomorphisms,}
% Trans. Amer. Math. Soc. {\bf 347} (1995) 1779-1786.
%%\bibitem[D2]{d2}
%%V. Deaconu, {\em A path model for circle algebras,} J. Operator
%%Theory {\bf 34} (1995) 57-89.
% \bibitem[D3]{d3}
% V. Deaconu, {\em Generalized Cuntz-Krieger algebras,}
% Proc. Amer. Math. Soc. {\bf 124} (1996) 3427-3435.
%%\bibitem[D4]{d4}
%%V. Deaconu, {\em Generalized solenoids and $C^*$-algebras,} in preparation.
\bibitem[Ki]{ki}
E. Kirchberg, {\em The classification of purely infinite $C^*$-algebras
using Kasparov's theory}, Fields Institute Communications, to appear.
%%Theorem C
\bibitem[Ku1]{ku1}
A. Kumjian, {\em  On $C^*$-diagonals}, Can. J. Math. {\bf 38} (1986)
969--1008.  % Imprimitivity groupoids considered in 5.7
\bibitem[Ku2]{ku2}
A. Kumjian, {\em  Notes on $C^*$-algebras of graphs}, in Operator
algebras and operator theory (Shanghai, 1997) Liming Ge et al., Eds.,
189--200, Contemp. Math., {\bf 228} Amer. Math. Soc., Providence, RI,
1998.
\bibitem[KMRW]{kmrw}
A. Kumjian, P. Muhly, J. Renault and D. Williams,
{\em The Brauer Group of a Locally Compact Groupoid},
       Amer. J. Math. {\bf 120} (1998) 901-954.
\bibitem[KP]{kp}
A. Kumjian and D. Pask, {\em Higher Rank Graph $C^*$-algebras},
New York Journal of Mathematics,  {\bf 6} (2000) 1-20.
% \bibitem[KPR]{kpr}
% A. Kumjian, D. Pask and I. Raeburn, {\em Cuntz-Krieger algebras
% of directed graphs,} Pacific J. Math., {\bf 184} (1998) 161-174.

\bibitem[M]{m}
P. Muhly, {\em Coordinates in Operator Algebra}.
CMBS Lecture Notes Series. (In progress).
% Theorem 5.31 Right Prinicpal H-space

\bibitem[MRW]{mrw} P. Muhly, J. Renault, and D. Williams,
{\em Equivalence and isomorphism for groupoid
$C^*$-algebras,} J. Operator Theory {\bf 17} (1987) 3-22.

% \bibitem[MS1]{ms1}
% P. Muhly and B. Solel, {\em On the simplicity of some Cuntz-Pimsner
% algebras,} Math. Scand., in press.
\bibitem[MW1]{mw1}
P. Muhly and D. Williams, {\em Continuous trace groupoid
$C^*$-algebras,} Math. Scand. {\bf 66} (1990) 231-241.
% Proposition 2.2 Proper principal groupoids
\bibitem[MW2]{mw2}
P. Muhly and D. Williams, {\em Groupoid cohomology and the
Dixmier-Douady class,} Proc. London Math.  Soc {\bf 71} (1995)
109-134.
% Theorem 3.5 Right Principal H-space bb
\bibitem[OPT]{opt}
D. Olesen, G. Pedersen, M. Takesaki,
{\em Ergodic actions of compact abelian groups}, J. Operator Theory
{\bf 3} (1980) 237--269.
\bibitem[Pd]{pd}
G. K. Pedersen, {\em $C^*$-algebras and their automorphism groups,}
Academic Press, London, 1979.
\bibitem[Ph]{ph}
N. C. Phillips, {\em A classification theorem for nuclear purely
infinite simple $C^*$-algebras}, Doc. Math. \textbf{5} (2000) 49--114.
%% Cor. 4.2.2

\bibitem[Pt1]{pt1}
I. Putnam, {\em $C^*$-algebras from Smale spaces}, Canad. J. Math. {\bf
48} (1996) 175-195.
\bibitem[Pt2]{pt2}
I. Putnam, {\em Hyperbolic systems and generalized Cuntz-Krieger
algebras,} Lecture notes from the summer school in operator algebras
(Odense, 1996).
\bibitem[PtS]{pts}
I. Putnam and J. Spielberg, {\em On the structure of $C^*$algebras
associated to hyperbolic dynamical systems,} J. Funct. Anal. {\bf 163}
(1999) 279-299.
\bibitem[Rn1]{rn1}
J. Renault, {\em A Groupoid approach to $C^*$-algebras,} Lecture Notes in
Mathematics {\bf 793}, Springer-Verlag, Berlin, 1980.
\bibitem[Rn2]{rn2}
J. Renault, {\em Repr\'{e}sentation des produits crois\'{e}s
d'alg\`{e}bres de groupo\"{\i}des,} J. Operator Theory {\bf 18} (1987)
67-97.
% pi systems - \i is ``dot-less i''
\bibitem[RSc]{rsc} J. Rosenberg and C. Schochet, {\em The K\"{u}nneth
theorem and the universal coefficient theorem for Kasparov's
generalized
$K$-functor,} Duke Math. J. {\bf 55} (1987), 431-474.
\bibitem[RSt1]{rst1}
G. Robertson and T. Steger, {\em $C^*$-algebras arising from group
actions on the boundary of a triangle building,} Proc. London
Math. Soc. {\bf 72} (1996) 613-637.
\bibitem[RSt2]{rst2}
G. Robertson and T. Steger, {\em Affine buildings, tiling systems and
higher rank Cuntz-Krieger algebras,} J. Reine Angew. Math., {\bf 513}
(1999) 115-144.
\bibitem[Ru1]{ru1}
D. Ruelle, {\em Thermodynamic Formalism}, Encyclopedia of Math. and
its appl. Vol. {\bf 5}, Addison-Wesley, Reading, 1978.
\bibitem[Ru2]{ru2}
D. Ruelle, {\em Noncommutative algebras for hyperbolic
diffeomorphisms}, Invent. Math. {\bf 93} (1988) 1-13.
\bibitem[Se]{sn}
E. Seneta, {\em Non-negative Matrices and Markov Chains}, Springer
Series in Statistics ($2^{\text{nd}}$ edition), Springer, New York,
1981.
\bibitem[Sm]{sm}
S. Smale, {\em Differentiable dynamical systems}, Bull. Amer. Math
Soc. {\bf 73} (1967) 747-817.
\bibitem[Z]{z}
S. Zhang, {\em Certain $C^*$-algebras with real rank zero and their
corona and multiplier algebras, Part I}, Pacific J. Math. {\bf 155}
(1992), 169-197. %% Theorem 1.2 Zhang's dichotomy
\end{thebibliography}
\end{document}